\documentclass[12pt]{amsart}

\usepackage{multirow}%
\usepackage{amsfonts}%
\usepackage{mathrsfs}%
\usepackage[title]{appendix}%
\usepackage{textcomp}%
\usepackage{manyfoot}%
\usepackage{algorithm}%
\usepackage{algorithmicx}%
\usepackage{algpseudocode}%
\usepackage{listings}%

\usepackage{scalerel}
\usepackage{mathtools}

\theoremstyle{plain}
\newtheorem{theorem}{Theorem}[section]

\theoremstyle{definition}
\newtheorem{definition}[theorem]{Definition}

\theoremstyle{remark}
\newtheorem{remark} [theorem]{Remark}

\numberwithin{equation}{section}

\makeatletter
\@namedef{subjclassname@2020}{%
	\textup{2020} Mathematics Subject Classification}
\makeatother

\begin{document}

%------
% Insert the title of your paper and (if necessary)
% a short title for the running head.
%------
\title{Lipschitz vector spaces}

\author[T. Valent]{Tullio Valent}
\address{Department of Mathematics ``Tullio Levi Civita''\\ University of Padua\\
	Via Trieste 63\\
	35121 Padua, Italy}
\email{tullio.valent@unipd.it}

\date{}

%------
% Add MSC 2020 codes according to www.ams.org/msc/msc2020.html.
% Secondary codes (in square brackets) are optional.
%------

%------
% Add a list of keywords. Only capitalise those keywords that start with a proper name.
%------

%------
% Insert your abstract.
%------
\begin{abstract}
The initial part of this paper is devoted to the notion of pseudo-seminorm on a vector space $E$. We prove that the topology of every topological vector space is defined by a family of pseudo-seminorms (and so, as it is known, it is uniformizable). Then we devote ourselves to the  Lipschitz vector structures on $E$, that is those Lipschitz structures on $E$ for which the addition is a Lipschitz map, while the scalar multiplication is a locally Lipschitz map, and we prove that any topological vector structure on $E$ is associated to some  Lipschitz vector structure.

Afterwards, we attend to the bornological Lipschitz maps. The final part of the article is devoted to the  Lipschitz vector structures compatible with locally convex topologies on $E$.
\end{abstract}



\maketitle

%------
% INSERT THE BODY OF THE PAPER HERE (except
% acknowledgments, funding info and bibliography)
%------
\section{Introduction}\label{sec1}

Sections \ref{sec2} and \ref{sec4} are the hard part of this paper. In Sect. \ref{sec2} we give the notion of pseudo-seminorm on the vector spaces, and we prove that the topology of every topological vector space is defined by a family of pseudo-seminorms, and so, among other things, it is uniformizable (as it is known).

A result (only) formally analogous to this one holds in the context of the topological spaces, because every topology is defined by a family of weak pseudo-metrics (see \cite{bib3}).
Sect. \ref{sec4} is devoted to the  Lipschitz vector structures on a vector space $E$, i.e., to the Lipschitz structures on $E$ such that the addition is a Lipschitz map, while the scalar multiplication is a locally Lipschitz map. We show that any topological vector structure on $E$ is (canonically) associated to some  Lipschitz vector structure, and vice versa. Sect. \ref{sec5} attends to the meaning of bornological Lipschitz maps. Finally, Sect. \ref{sec6} is devoted to the  Lipschitz vector structures compatible with locally convex topologies on $E$. 

\section{Pseudo-Seminorms}\label{sec2}

Let $E$ be a vector space over the scalar field $\mathbb{K}$ (with $\mathbb{K}=\mathbb{R}$ or $\mathbb{K}=\mathbb{C}$).
\begin{definition}\label{pseudoseminorm}
	Any map $\mathbf{x} \mapsto |\mathbf{x}|$ from $E$ to $\mathbb{R}$ such that 
	\begin{equation}\label{eq:pseudoseminorm}
		\begin{cases*}
			|\lambda \mathbf{x}| \leq |\mathbf{x}| \quad \text{for all} \quad \mathbf{x} \in E \quad \text{and} \quad |\lambda| \leq 1\\
			|\mathbf{x}_1 + \mathbf{x}_2| \leq |\mathbf{x}_1| + |\mathbf{x}_2| \quad \text{for all} \quad \mathbf{x}_1,\mathbf{x}_2 \in E\\
			\vert\mathbf{0} \vert = 0
		\end{cases*}
	\end{equation} 
	(where $\mathbf{0}$ denotes the origin of $E$) will be called a \emph{pseudo-seminorm} on $E$. 
\end{definition}
Observe that $|\mathbf{x}| \geq 0$ for all $\mathbf{x} \in E$, because \eqref{eq:pseudoseminorm} implies $0=|0|=|0\mathbf{x}| \leq |\mathbf{x}|$.

It is easy to check that the map $(\mathbf{x}_1,\mathbf{x}_2) \mapsto |\mathbf{x}_1-\mathbf{x}_2|$, from $E \times E$ to $\mathbb{R}$, is a pseudo-metric for $E$ that is invariant under translations.

\begin{theorem}\label{th:pseudoseminorm}
	If $(V_n)_{n \in \mathbb{N}}$ is a sequence of circled subsets of $E$ such that
	\begin{equation}\label{eq:circ}
		V_{n+1} + V_{n+1} \subseteq V_n \quad \text{for all} \quad n \in \mathbb{N},
	\end{equation}
	then there exists a pseudo-seminorm $\mathbf{x} \mapsto |\mathbf{x}|$ on $E$ such that 
	\begin{equation}\label{eq:pseudocirc}
		\left\{\mathbf{x} \in E : |\mathbf{x}|< \dfrac{1}{2^n}\right\} \subseteq V_n \subseteq \left\{\mathbf{x} \in E : |\mathbf{x}|\leq \dfrac{1}{2^n}\right\} \quad
		\text{for all} \quad n \in \mathbb{N}.
	\end{equation} 
\end{theorem}

\begin{proof}
	For each finite and nonempty subset $I$ of $\mathbb{N}$ we set
	\begin{equation*}
		V_I = \sum_{i \in I} V_i
	\end{equation*}
	and we define the function $\mathbf{x} \mapsto |\mathbf{x}|$ from $E$ to $[0,1]$ by setting
	\begin{equation*}
		\begin{cases*}
			|\mathbf{x}| = 1 \quad \text{if} \quad x \not\in V_I \quad \text{for every} \quad V_I\\
			|\mathbf{x}| = \inf\left\{\sum\limits_{i \in I}\dfrac{1}{2^i}: x \in V_I \right\} \quad \text{if} \quad x \in V_I \quad \text{for some} \quad I.
		\end{cases*}
	\end{equation*}
	
	Let us now to prove that this function is a pseudo-seminorm. The first property of \eqref{eq:pseudoseminorm} is satisfied because every $V_I$ is circled, while the third property of \eqref{eq:pseudoseminorm} is satisfied because $\mathbf{0} \in V_n$ for all $n \in \mathbb{N}$. Then now we attend to the second condition of \eqref{eq:pseudoseminorm}. Put $p_I= \sum_{i \in I} \tfrac{1}{2^i}$. As the second condition is obviously satisfied when $|\mathbf{x}_1| + |\mathbf{x}_2| \geq 1$, we consider the case $|\mathbf{x}_1| + |\mathbf{x}_2| < 1$. If $\varepsilon$ is a real number $>0$ such that $|\mathbf{x}_1| + |\mathbf{x}_2| +2\varepsilon< 1$, there are two finite and non-empty subsets $I_1$ and $I_2$ of $\mathbb{N}$ such that $\mathbf{x}_1 \in V_{I_1}$, $\mathbf{x}_2 \in V_{I_2}$, and $p_{I_1} < |\mathbf{x}_1| + \varepsilon$, $p_{I_2} < |\mathbf{x}_2| + \varepsilon$.
	
	Since  $p_{I_1} + p_{I_2} < 1$ there is an unique and non-empty subset $I$ of $\mathbb{N}$ such that $p_{I_1} + p_{I_2} = p_I$. Hence \eqref{eq:circ} implies 
	$ V_{I_1}+ V_{I_2} \subseteq V_I$. Then $\mathbf{x}_1 + \mathbf{x}_2 \in V_I$. Thus $|\mathbf{x}_1 + \mathbf{x}_2| \leq p_I = p_{I_1} + p_{I_2} < |\mathbf{x}_1| + |\mathbf{x}_2|+ 2\varepsilon$, and so (being $\varepsilon$ arbitrarily small) $|\mathbf{x}_1 + \mathbf{x}_2| \leq |\mathbf{x}_1| + |\mathbf{x}_2|$. 
	
	Nothing remains but to show that \eqref{eq:pseudocirc} is satisfied. If $\mathbf{x} \in V_n$, then (evidently) $|\mathbf{x}| \leq \tfrac{1}{2^n}$. On the other hand, if $|\mathbf{x}| <   \tfrac{1}{2^n}$ there is a finite and non-empty subset $I$ of $\mathbb{N}$ such that  $\mathbf{x} \in V_I$ and $\sum_{i \in I} \tfrac{1}{2^i} < \tfrac{1}{2^n}$. This implies $n < i$ for all $i \in I$, which implies $V_i \subseteq V_n$ for all $i \in I$, namely $V_I \subseteq V_n$, and so $\mathbf{x} \in V_n$.
\end{proof}

\begin{theorem}\label{th:seminorm}
	The topology of every topological vector space $E$ is defined by a family of pseudo-seminorms (i.e., by the family of pseudo-metrics associated to such pseudo-seminorms), and so it is uniformizable.
\end{theorem}
\begin{proof}
	Let $\mathscr{F}(\mathbf{0})$ denote the filter of neighborhoods of $\mathbf{0}$ in $E$. For every $V \in \mathscr{F}(\mathbf{0})$ choose a sequence $(V_n)_{n \in \mathbb{N}}$ of circled elements of $\mathscr{F}(\mathbf{0})$ such that
	\begin{equation*}
		V_1 \subseteq V, \, \, V_{n+1} +V_{n+1} \subseteq V_n \quad \text{for} \quad n >1.
	\end{equation*}
	It is well known that such a sequence exists in any topological vector space. By Theorem~\ref{th:pseudoseminorm} there is a pseudo-seminorm on $E$ that we will simply denote by $|\cdot|_{\scriptscriptstyle V}$ (even if it depends, besides on $V$, also on the sequence $(V_n)_{n \in \mathbb{N}}$) such that
	\begin{equation}\label{eq:sequence}
		\left\{\mathbf{x} \in E : |\mathbf{x}|_{\scriptscriptstyle V} < \dfrac{1}{2^n}\right\} \subseteq V_n \subseteq 	\left\{\mathbf{x} \in E : |\mathbf{x}|_{\scriptscriptstyle V} \leq \dfrac{1}{2^n}\right\}.
	\end{equation}
	
	Let us denote by $\mathscr{T}'_{E}$ the topology on $E$ defined by the family $(|\cdot|_{\scriptscriptstyle V})_{V \in \mathscr{F}(\mathbf{0})}$, i.e. the topology on $E$ defined by the family $(d_{\scriptscriptstyle V})_{V \in \mathscr{F}(\mathbf{0})}$ of the translation invariant pseudo-metrics  $d_{\scriptscriptstyle V} : (\mathbf{x}_1, \mathbf{x}_2) \mapsto |\mathbf{x}_1 - \mathbf{x}_2|_{\scriptscriptstyle V}$. 
	
	A base of neighborhoods of $\mathbf{0}$ for $\mathscr{T}'_{E}$ is
	\begin{align*}
		&\left(\left\{\mathbf{x} \in E : |\mathbf{x}|_{\scriptscriptstyle V} < \dfrac{1}{2^n} \right\}, n \in \mathbb{N}, V \in \mathscr{F}(\mathbf{0})\right), \quad
		\text{or} \\
		&\left(\left\{\mathbf{x} \in E : |\mathbf{x}|_{\scriptscriptstyle V} \leq \dfrac{1}{2^n} \right\}, n \in \mathbb{N}, V \in \mathscr{F}(\mathbf{0})\right).
	\end{align*}
	From these we can obtain, by translations, two bases of neighborhoods of any point $\mathbf{x}$ of $E$. Then, in view of \eqref{eq:sequence}, the filter of the neighborhoods of $\mathbf{0}$ (and hence of every point of $E$) is the same for the topologies  $\mathscr{T}'_{E}$ and  $\mathscr{T}_{E}$, where  $\mathscr{T}_{E}$ denotes the topology of the topological vector space $E$.
	
	Then,  $\mathscr{T}'_{E}=  \mathscr{T}_{E}$, and thus the topology of any topological vector space $E$ can be defined by the family $\{ |\cdot|_{\scriptscriptstyle V}: V \in \mathscr{F}(\mathbf{0})\}$ of pseudo-seminorms.
	
\end{proof}

\section{Preliminaries about Lipschitz structures for any set}\label{sec3}
Before talking about Lipschitz structures for vector spaces, let us shortly remind some facts concerning the Lipschitz structures for any set. An axiomatic definition of Lipschitz structures for a set can be found in \cite{bib1}.

By a \emph{Lipschitz structure} for a set $X$ we mean a non-empty family $\mathscr{L}$ of pseudo-metrics $d$  on $X$ satisfying the following conditions:

\begin{itemize}
	\item[$(L_1)$] \quad $d \leq d_1, d_1 \in \mathscr{L} \Rightarrow d \in  \mathscr{L}$;
	\item[$(L_2)$] \quad  $d \in \mathscr{L} \Rightarrow \alpha d \in  \mathscr{L}$ for every real number $\alpha >0$;
	\item[$(L_3)$] \quad $d_1,d_2 \in \mathscr{L} \Rightarrow d_1 \vee d_2 \in \mathscr{L}$.
\end{itemize}

It is easy to prove that	$\left((L_1),(L_2),(L_3)\right)$ is equivalent to 	$\left((L_1),(L_4)\right)$, with \\
$(L_4) \quad d_1,d_2 \in \mathscr{L} \Rightarrow d_1 + d_2 \in \mathscr{L}$.

The pair $(X,\mathscr{L})$ will be called a \emph{Lipschitz space}. A base of  $\mathscr{L}$ is a subset $\mathscr{B}$ of $\mathscr{L}$ such that for every $d \in \mathscr{L}$ there are $b \in \mathscr{B}$ and $\alpha >0$ such that
$d \leq \alpha b.$
If $\mathscr{B}$ is a base of $\mathscr{L}$ then 
$$\mathscr{L} = \{d: d \leq \alpha b \quad \text{for some} \quad b \in \mathscr{B} \quad \text{and}\quad \alpha >0\}.$$

A set $\mathscr{B}$ of pseudo-metrics on $X$ is said a base for a Lipschitz structure for $X$ if there is  Lipschitz structure of which $\mathscr{B}$ is a base. This is true if and only if the set $\{d: d \leq \alpha b \quad \text{with} \quad b \in \mathscr{B} \quad \text{and}\quad \alpha >0\}$ is a Lipschitz structure for $X$.

It follows that $\mathscr{B}$ is a base for a Lipschitz structure for $X$ if and only if  for every $b_1,b_2 \in \mathscr{B}$ there are $b \in \mathscr{B}$ and $\alpha >0$ such that $b_1 \vee b_2 \leq \alpha b$. Consequently, given a non-empty family $\mathscr{P}$ of pseudo-metrics on $X$, the family $\mathscr{B}$ of the suprema of all finite subsets of $\mathscr{P}$ is a base for a Lipschitz structure $\mathscr{L}(\mathscr{P})$ for $X$. Of course, $\mathscr{L}(\mathscr{P})$ is the smallest  Lipschitz structure for $X$ containing $\mathscr{P}$, and it is called the \emph{Lipschitz structure for X generated by $\mathscr{P}$}. We have
$$\mathscr{L}(\mathscr{P})=\{d \in \mathscr{P}(X): d \leq \alpha (d_1\vee  \dots \vee d_n), d_1,\dots,d_n \in \mathscr{P}, n \geq 1, \alpha >0 \}.$$  	
The Lipschitz structure implicitly considered on $\mathbb{R}$, or $\mathbb{C}$, will be the one generated by the usual metric.

Let $(X,\mathscr{L}_X)$ and $(Y,\mathscr{L}_Y)$ be  Lipschitz spaces. A map $f: X \mapsto Y$ is called a \emph{Lipschitz map} if for every $d_Y \in \mathscr{L}_Y$ there is $d_X \in \mathscr{L}_X$ such that

$$d_Y(f(x_1),f(x_2))\leq d_X(x_1,x_2) \quad \forall x_1,x_2 \in X.$$

If $\mathscr{L}_Y$ is generated by  $\mathscr{P}_Y$, then a map $f:X \mapsto Y$ is a Lipschitz map whenever for every $d \in \mathscr{P}_Y$ there is $d_X \in \mathscr{L}_X$ such that 

$$d(f(x_1),f(x_2))\leq d_X(x_1,x_2) \quad \forall x_1,x_2 \in X.$$

\emph{To each Lipschitz structure $\mathscr{L}$ can be (canonically) associated a topological structure $\tau(\mathscr{L})$}. $\tau(\mathscr{L})$ is the topology of the uniformity generated by the family $\mathscr{L}$ of pseudo-metrics.

Let $(X,\mathscr{L}_X)$ and $(Y,\mathscr{L}_Y)$	be Lipschitz spaces. A map $f: X \mapsto Y$ will be called a \emph{locally Lipschitz map}  if every $x \in X$ has a neighborhood $U_x$ for the topological structure (canonically) associated to $\mathscr{L}_X$ such that for every $d_Y \in \mathscr{L}_Y$ there is $d_X \in \mathscr{L}_X$ such that

\begin{equation*}\label{funceq2}
	d_Y(f(x_1),f(x_2)) \leq d_X(x_1,x_2) \quad \forall x_1,x_2 \in U_x.
\end{equation*} 

Let again $(X,\mathscr{L}_X)$ and $(Y,\mathscr{L}_Y)$	be Lipschitz spaces. The \emph{product Lipschitz structure}, $\mathscr{L}_{X \times Y}$, for $X \times Y$ is the smallest Lipschitz structure for which the projections $(x,y) \mapsto x$ and $(x,y) \mapsto y$ (from $X \times Y$ to $X$, and from $X \times Y$ to $Y$) are Lipschitz maps. $\mathscr{L}_{X \times Y}$ is generated by the family of the pseudo-metrics $d$ on $X \times Y$ of the type
$$d((x_1,y_1),(x_2,y_2)) \coloneqq d_X(x_1,x_2) + d_Y(y_1,y_2),$$
with $d_X \in \mathscr{L}_X$ and $d_Y \in \mathscr{L}_Y$.

\section{Lipschitz vector structures}\label{sec4}

\begin{definition}
	The Lipschitz structure for a topological vector space $E$ generated by the family of the pseudo-metrics $d_U$, where $U$ is a neighborhood of $\mathbf{0}$ in $E$, and $d_U$ is obtained from the pseudo-seminorm $p_U$ considered in proof of Theorem~\ref{th:seminorm}, will be called the \emph{Lipschitz structure} (canonically) \emph{associated to the  topological vector structure of $E$}. 
\end{definition}
\begin{definition}
	A Lipschitz structure for a  vector space $E$ will be called a  \emph{Lipschitz vector structure} if the addition is a Lipschitz map and the scalar multiplication is a locally Lipschitz map. A vector space endowed with a  Lipschitz vector structure will be called a \emph{Lipschitz vector space}. 
\end{definition}

\begin{remark}
	The topological structure (canonically) associated to a  Lipschitz vector structure $\mathscr{L}$ for a vector space $E$ (in the sense that it is obtained from the uniformity defined by the family $\mathscr{L}$ of pseudo-metrics) is the one of a topological vector space, because each locally Lipschitz map is continuous (for the topological structures associated to the Lipschitz structures).  Thus \emph{every Lipschitz vector space can be} (canonically) \emph{associated to some topological vector space}. The following theorem asserts that vice versa \emph{every topological vector structure is} (canonically) \emph{associated to some  Lipschitz vector structure}.
\end{remark}

\begin{theorem}
	Let $(E,\tau)$ be any topological vector space. The Lipschitz structure $\mathscr{L}(\tau)$ for $E$ (canonically) associated to $\tau$ is a  Lipschitz vector structure.
\end{theorem}
\begin{proof}
	The Lipschitz structure $\mathscr{L}(\tau)$ is generated by the family $\{d_U: U \in \mathscr{F}(\mathbf{0})\}$, where $\mathscr{F}(\mathbf{0})$ denotes the filter of the neighborhoods of the origin $\mathbf{0}$ for the topological vector space $(E,\tau)$, and $d_U$ is the pseudo-metric defined by $d_U(\mathbf{x}_1,\mathbf{x}_2) \coloneqq p_U(\mathbf{x}_1 - \mathbf{x}_2)$ for all $\mathbf{x}_1, \mathbf{x}_2 \in E$, with $p_U$ the pseudo-seminorm that, in the proof of Theorem~\ref{th:seminorm} is denoted by $|\cdot|_{\scriptscriptstyle U}$. We start by remarking that the Lipschitz structure $\mathscr{L}(\tau)$ is the family of all pseudo-metrics $d$ for $E$ such that $d \leq d_U$ for some $U \in \mathscr{F}(\mathbf{0})$. This is true because $d_{U_1} \vee d_{U_2} \leq d_{U_1 \cup U_2}$ for all $U_1,U_2 \in \mathscr{F}(\mathbf{0})$, and $|\lambda|d_U = d_{|\lambda|^{-1}V}$ for all scalar $\lambda \not=0$, with $V$ a circled element of $\mathscr{F}(\mathbf{0})$ contained in $U$.
	
	It follows that the topological structure on the vector space $E$ associated to the Lipschitz structure generated by the family $\{d_U : U \in \mathscr{F}(\mathbf{0})\}$ of pseudo-metrics is defined by the family $\{p: p \leq p_U \, \, \text{for some} \, \, U \in \mathscr{F}(\mathbf{0})\}$ of pseudo-seminorms $p$ on $E$.
	
	Then, as $\{\mathbf{x} \in E : p_U(\mathbf{x}) \leq \varepsilon\} \subseteq \{\mathbf{x} \in E : p(\mathbf{x}) \leq \varepsilon\}$ for all $U \in \mathscr{F}(\mathbf{0})$ and $\varepsilon$ a real number $>0$, the filter of the neighborhoods of $\mathbf{0}$ for the topological structure on $E$ associated to the Lipschitz structure on $E$ generated by the family $\{d_U : U \in \mathscr{F}(\mathbf{0})\}$ of pseudo-metrics is the same as $\mathscr{F}(\mathbf{0})$.
	
	We now attend to the proof of the theorem, starting by proving that the addition is a Lipschitz map (from $E \times E$ to $E$) when the Lipschitz structure on $E$ is generated by the family $\{d_U : U \in \mathscr{F}(\mathbf{0})\}$ of pseudo-metrics. We must observe that the product Lipschitz structure for $E \times E$ is generated by the family of the pseudo-metrics $d$ on $E \times E$ of the type
	\begin{equation}\label{eq:lipscprod}
		d((\mathbf{x}_1,\mathbf{y}_1),(\mathbf{x}_2,\mathbf{y}_2)) \coloneqq d_{U_1}(\mathbf{x}_1,\mathbf{x}_2) + d_{U_2}(\mathbf{y}_1,\mathbf{y}_2)
	\end{equation}
	with $U_1,U_2 \in \mathscr{F}(\mathbf{0})$. Therefore, the addition $\left((\mathbf{x}_1,\mathbf{y}_1),(\mathbf{x}_2,\mathbf{y}_2)\right) \mapsto (\mathbf{x}_1+\mathbf{x}_2,\mathbf{y}_1+\mathbf{y}_2)$ is a Lipschitz map if for every $U \in \mathscr{F}(\mathbf{0})$ there are $U_1,U_2 \in \mathscr{F}(\mathbf{0})$ such that
	\begin{equation*}
		d_U(\mathbf{x}_1+\mathbf{x}_2,\mathbf{y}_1+\mathbf{y}_2) \leq d_{U_1}(\mathbf{x}_1,\mathbf{x}_2) + d_{U_2}(\mathbf{y}_1,\mathbf{y}_2)
	\end{equation*} 
	for all $\mathbf{x}_1,\mathbf{x}_2,\mathbf{y}_1,\mathbf{y}_2 \in E$. In our case it suffices to take $U_1=U_2=U$, for
	
	\begin{equation*}
		d_U(\mathbf{x}_1+\mathbf{x}_2,\mathbf{y}_1+\mathbf{y}_2) = d_U((\mathbf{x}_1,\mathbf{y}_1)+(\mathbf{x}_2,\mathbf{y}_2)) \leq
		d_U(\mathbf{x}_1,\mathbf{y}_1) +  d_U(\mathbf{x}_2,\mathbf{y}_2).
	\end{equation*}
	Hence we can conclude that the addition is a Lipschitz map. 
	
	Now we observe that from the arguments developed in Section~\ref{sec2} it is possible to show that if $U,V \in \mathscr{F}(\mathbf{0})$ and $U=sV$, with $s$ a positive number, then
	\begin{equation}\label{eq:scalar}
		|\lambda| \leq s \Rightarrow p_U(\lambda \mathbf{x}) \leq p_V(\mathbf{x}) \quad \forall x \in E. 
	\end{equation}
	We can deduce that for every $(\bar{\lambda},\bar{\mathbf{x}}) \in \mathbb{K} \times E$ and every $U \in \mathscr{F}(\mathbf{0})$ there are $V \in \mathscr{F}(\mathbf{0})$ and a number $r >0$ such that
	\begin{equation*}
		d_U(\lambda_1 \mathbf{x}_1,\lambda_2 \mathbf{x}_2) \leq d_V(\mathbf{x}_1,\mathbf{x}_2)
	\end{equation*}
	whenever $|\lambda_1 -\bar{\lambda}| \leq r, |\lambda_2 -\bar{\lambda}| \leq r$. Indeed, putting $s=r + |\bar{\lambda}|$, the inequality $|\lambda_1-\bar{\lambda}| \leq r$ implies $|\lambda_1| \leq s$. Then, by \eqref{eq:scalar}
	
	\begin{align*}
		d_U(\lambda_1,\mathbf{x}_1,\lambda_2 \mathbf{x}_2)& \leq d_U(\lambda_1\mathbf{x}_1,\lambda_1 \mathbf{x}_2)+ d_U(\lambda_1\mathbf{x}_2,\lambda_2 \mathbf{x}_2)=\\
		&= p_U(\lambda_1(\mathbf{x}_1-\mathbf{x}_2))+p_U((\lambda_1+ \lambda_2)(\mathbf{x}_2-\mathbf{x}_2))=\\ 
		&= p_U(\lambda_1(\mathbf{x}_1-\mathbf{x}_2)) \leq p_V(\mathbf{x}_1-\mathbf{x}_2)) =d_V(\mathbf{x}_1,\mathbf{x}_2).
	\end{align*}
	
	Thus  we can clearly conclude that the scalar multiplication $(\lambda,\mathbf{x}) \mapsto \lambda \mathbf{x}$ is a locally Lipschitz map. Hence the Lipschitz structure on $E$ generated by the family $\{d_U: U \in \mathscr{F}(\mathbf{0})\}$ of pseudo-metrics is a  Lipschitz vector structure.
\end{proof}

\begin{definition}
	Let $E,F$ be topological vector spaces. A map $f : E \mapsto F$ will be called a \emph{Lipschitz map} if it is a Lipschitz map for the  Lipschitz vector structures (canonically) associated to the structure of topological vector space on $E$ and $F$.
\end{definition}
Thus, $f$ is a Lipschitz map if and only if for every $U \in \mathscr{F}_F(\mathbf{0})$ there are $V \in \mathscr{F}_E(\mathbf{0})$ and a number $c >0$ such that
\begin{equation*}
	p_U(f(\mathbf{x}_1)-f(\mathbf{x}_2)) \leq c p_V(\mathbf{x}_1-\mathbf{x}_2) \quad \text{for all}\quad \mathbf{x}_1,\mathbf{x}_2 \in E,
\end{equation*} 
where $\mathscr{F}_E(\mathbf{0})$ and $\mathscr{F}_F(\mathbf{0})$ denote the filters of the neighborhoods of $\mathbf{0}$ in $E$ and $F$, and $p_U,p_V$ are the pseudo-seminorms defined in Section~\ref{sec2}.

\section{Bornological Lipschitz maps}\label{sec5}
Let $(E, \mathscr{L}_E)$ be a Lipschitz vector space, and $A$ be a convex, circled and bounded subset of $E$. On the smallest vector subspace $E_A$ of $E$ containing $A$ consider the seminorm $p_A$ (the Minkowski functional, see \cite[p.~15]{bib2}) defined by
\begin{equation*}
	p_A(\mathbf{x}) \coloneqq \inf \{\lambda >0: \mathbf{x} \in \lambda A\}, \quad \mathbf{x} \in E_A.
\end{equation*}
The Lipschitz structure on $E_A$ generated by the pseudo-metric $d_A$ associated to $p_A$ will be denoted by $\mathscr{L}(p_A)$.

\begin{theorem}
	If the Lipschitz structure $\mathscr{L}_E$ is generated by a family $\mathscr{P}$ of seminorms (so giving a locally convex topology on $E$), then $\mathscr{L}_E$ induces on $E_A$ a Lipschitz structure contained in   $\mathscr{L}(p_A)$.
\end{theorem}
\begin{proof}
	We first observe that the elements of $E_A$ are of type $\lambda \mathbf{a}$ with $\lambda \geq 0$ and $\mathbf{a} \in A$. We must prove that for every $p \in \mathscr{P}$ there is a number $l >0$ such that $p(\mathbf{x}) \leq l p_A(\mathbf{x})$ for all $\mathbf{x} \in E_A$. This is true with $l=\sup p(A)$ [which is a real number, because $A$ is bounded]. Indeed, for every $\mathbf{x}$ ($=\lambda \mathbf{a}$, with $\lambda \geq 0$ and $\mathbf{a} \in A$) of $E_A$, we have
	\begin{equation*}
		p(\mathbf{x}) = \lambda p(\mathbf{a}) \leq l \lambda \leq l p_A(\mathbf{x}),
	\end{equation*} 
	and so, if $d$ denotes the pseudo-metric associated to $p$,
	\begin{equation*}
		d(\mathbf{x}_1,\mathbf{x}_2) \leq l p_A(\mathbf{x}_1-\mathbf{x}_2) \quad \text{for all} \quad \mathbf{x}_1,\mathbf{x}_2 \in E_A.
	\end{equation*}
\end{proof}

\begin{definition}
	A map $f : E \mapsto X$, with $E$ a Lipschitz vector space and $X$ a Lipschitz space, will be said a \emph{bornological Lipschitz map} if for every convex, circled and bounded subset $A$ of $E$, $f$ is a Lipschitz map from $E_A$, endowed with the Lipschitz structure generated by the pseudo-metric associated to $p_A$, to $X$.
\end{definition}

\begin{theorem}
	Let $(E, \mathscr{L}_E)$ be a Lipschitz vector space with the Lipschitz structure $\mathscr{L}_E$ generated by a family $\mathscr{P}$ of seminorms, and let $(X,\mathscr{L}_X)$ be a Lipschitz space. Every Lipschitz map $f:E \mapsto X$ is a bornological Lipschitz map.
\end{theorem}

\begin{proof}
	This theorem in an immediate consequence of the previous one, because (by virtue of the previous Theorem) for every convex, circled and bounded subset $A$ of $E$, the Lipschitz structure of the subspace $E_A$ of $E$ is contained in the Lipschitz structure generated by the pseudo-metric associated to $p_A$.
	
\end{proof}

\section{Lipschitz vector structures compatible with locally convex topologies}\label{sec6}

In the case of locally convex topologies the treatment is, evidently, easier than in the case of a general vector topology. Moreover, the  Lipschitz vector structures are generated by families of pseudo-metrics associated to seminorms (not only to pseudo-seminorms, as it occurs in the general case).

\begin{theorem}
	Let $(E,\tau)$ be a locally convex topological vector space, and let $\mathscr{L}(\tau)$ be the Lipschitz structure (canonically associated to $\tau$) generated by the family $\{d_U: U \text{convex neighborhood of} \,\, \mathbf{0} \,\, \text{for} \, \, \tau\}$, where $d_U$ is the pseudo-metric defined by $d_U(\mathbf{x}_1,\mathbf{x}_2)=p_U(\mathbf{x}_1-\mathbf{x}_2)$, with $p_U$ the seminorm defined by
	\begin{equation}\label{eq:Minkow}
		p_U(\mathbf{x}) \coloneqq \inf \{\lambda >0 : \mathbf{x} \in \lambda U\}.
	\end{equation}
	$\mathscr{L}(\tau)$ is a  Lipschitz vector structure on E [i.e., when $E$ is endowed with the Lipschitz structure $\mathscr{L}(\tau)$ the addition is a Lipschitz map from $E \times E$ to $E$, and the scalar multiplication from $\mathbb{K} \times E$ to $E$ is a locally Lipschitz map], and $\tau$ is the topology (canonically) associated to $\mathscr{L}(\tau)$.
\end{theorem}

\begin{proof}
	We first observe that, as $U$ is convex, the (Minkowski) functional) $p_U$ defined by \eqref{eq:Minkow} is a seminorm. We also remark that the Lipschitz structure  $\mathscr{L}(\tau)$ is the family of all pseudo-metrics $d$ for $E$ such that $d \leq d_U$ for some convex neighborhood $U$ of $\mathbf{0}$ for $\tau$. This implies, without difficulty, that the topology (canonically) associated to $\mathscr{L}(\tau)$ is $\tau$.
	
	Let us now to show that the addition from $E \times E$ to $E$ is a Lipschitz map. Since the product Lipschitz structure on $E \times E$ is generated by the family of the pseudo-metrics of the type \eqref{eq:lipscprod}, with $U_1,U_2$ convex neighborhoods of $\mathbf{0}$ in $E$, the addition is a Lipschitz map because, for every convex neighborhood $U$ of $\mathbf{0}$ in $E$ we have
	
	\begin{align*}
		&d_U(\mathbf{x}_1+\mathbf{x}_2,\mathbf{y}_1 +\mathbf{y}_2) = p_U((\mathbf{x}_1+\mathbf{x}_2)-(\mathbf{y}_1+\mathbf{y}_2))=\\
		&=p_U((\mathbf{x}_1-\mathbf{y}_1)+(\mathbf{x}_2-\mathbf{y}_2)) \leq p_U(\mathbf{x}_1-\mathbf{y}_1) + p_U(\mathbf{x}_2-\mathbf{y}_2)=\\
		&=d_U(\mathbf{x}_1,\mathbf{y}_1)+d_U(\mathbf{x}_2,\mathbf{y}_2).
	\end{align*}
	Finally, we will prove that the scalar multiplication from $\mathbb{K} \times E$ to $E$ is a locally Lipschitz map. Indeed, we will show that for every $(\bar{\lambda},\bar{\mathbf{x}}) \in \mathbb{K} \times E$ and every convex neighborhood $U$ of $\mathbf{0}$ in $E$ there are two positive numbers $c_1$ and $c_2$ such that
	\begin{equation}\label{eq:locconv1}
		d_U(\lambda_1 \mathbf{x}_1,\lambda_2 \mathbf{x}_2) \leq c_1 d_U(\mathbf{x}_1,\mathbf{x}_2)+c_2 |\lambda_1-\lambda_2|
	\end{equation} 
	whenever
	\begin{equation}\label{eq:locconv2}
		|\lambda_1-\bar{\lambda}| \leq 1, \, |\lambda_2-\bar{\lambda}| \leq 1, \, d_U(\mathbf{x}_1,\bar{\mathbf{x}}) \leq 1, \, d_U(\mathbf{x}_2,\bar{\mathbf{x}}) \leq 1.
	\end{equation}
	
	To prove this we observe that 
	\begin{align*}
		&d_U(\lambda_1 \mathbf{x}_1,\lambda_2 \mathbf{x}_2)=p_U(\lambda_1 \mathbf{x}_1-\lambda_2 \mathbf{x}_2) =p_U((\lambda_1-\lambda_2)\mathbf{x}_1+\lambda_2(\mathbf{x}_1-\mathbf{x}_2)) \leq\\
		&\leq p_U((\lambda_1-\lambda_2)\mathbf{x}_1) +p_U(\lambda_2(\mathbf{x}_1-\mathbf{x}_2)) \leq\\
		&\leq |\lambda_1-\lambda_2|p_U(\mathbf{x}_1)+|\lambda_2|p_U(\mathbf{x}_1-\mathbf{x}_2),
	\end{align*}
	
	and so, by \eqref{eq:locconv2} (which implies $p_U(\mathbf{x}_1) \leq p_U((\mathbf{x}_1-\bar{\mathbf{x}})+\bar{\mathbf{x}}) \leq p_U(\mathbf{x}_1-\bar{\mathbf{x}})+p_U(\bar{\mathbf{x}}) \leq 1 + p_U(\bar{\mathbf{x}})$) we obtain
	\begin{equation*}
		d_U(\lambda_1 \mathbf{x}_1,\lambda_2 \mathbf{x}_2) \leq |\lambda_1- \lambda_2| (1+p_U(\bar{\mathbf{x}})) + (1+|\bar{\lambda}|)d_U(\mathbf{x}_1,\mathbf{x}_2).
	\end{equation*}
	Therefore \eqref{eq:locconv1} is true with $c_1 = 1 + |\bar{\lambda}|$ and $c_2 = 1 +p_U(\bar{\mathbf{x}})$.
\end{proof}

\end{document}